\documentclass{article}
\usepackage{amsmath}
\usepackage{latexsym}
\usepackage{amssymb,amsthm,amsxtra}
\usepackage{multicol}

\def\[#1\]{\begin{equation}#1\end{equation}}
\makeatletter
\def\beq{%
   \relax\ifmmode
      \@badmath
   \else
      \ifvmode
         \nointerlineskip
         \makebox[.6\linewidth]%
      \fi
      $$%%$$ BRACE MATCH HACK
   \fi
}
\def\eeq{%
   \relax\ifmmode
      \ifinner
         \@badmath
      \else
         $$%%$$ BRACE MATCH HACK
      \fi
   \else
      \@badmath
   \fi
   \ignorespaces
}

\def\enddisplaymath{\eeq\global\@ignoretrue}
\makeatother

\newtheorem{thm}{Theorem}
\newtheorem{cor}[thm]{Corollary}
\newtheorem{lem}[thm]{Lemma}

\theoremstyle{remark}
\newtheorem*{rem}{Remark}

\theoremstyle{definition}

\numberwithin{thm}{section}
\renewcommand{\Re}{\operatorname{Re}}

\DeclareMathOperator{\eig}{eig}
\newcommand{\wt}{\widetilde}

\newcommand{\R}{\mathbb R}

\newcommand{\Z}{\mathbb Z}

\newcommand{\fa}{\forall}

\begin{document}

\title{Images of eigenvalue distributions under power maps}
\author{
Eric M. Rains\\AT\&T Research\\rains@research.att.com}
\date{August 10, 2000}
\maketitle

\begin{abstract}
In \cite{Rains:1997a}, it was shown that if $U$ is a random $n\times n$
unitary matrix, then for any $p\ge n$, the eigenvalues of $U^p$ are
i.i.d. uniform; similar results were also shown for general compact
Lie groups.  We study what happens when $p<n$ instead.  For the
classical groups, we find that we can describe the eigenvalue distribution
of $U^p$ in terms of the eigenvalue distributions of smaller classical
groups; the earlier result is then a special case.  The proofs rely
on the fact that a certain subgroup of the Weyl group is itself a Weyl
group.  We generalize this fact, and use it to study the power-map
problem for general compact Lie groups.
\end{abstract}

In \cite{Rains:1997a}, it was shown that if a (uniformly) random $n\times
n$ unitary matrix $U$ is raised to a power $p\ge n$, the eigenvalues of the
resulting matrix are (exactly) independently distributed; this despite the
rather complicated dependence between the eigenvalues of $U$ itself.
Our purpose in the present note is to extend this result to the case $p<n$.
We find that the eigenvalue distribution of $U^p$ can in that case be
described in terms of a union of $p$ independent distributions, each
of which is itself the eigenvalue distribution of a random unitary matrix.
More precisely, we have:
\[
U(n)^p \sim \bigoplus_{0\le i<p} U(\lceil\frac{n-i}{p}\rceil).
\]
That is, if we take the $p$th power of a uniformly distributed element of
$U(n)$, the resulting eigenvalue distribution is the same as if we
took the union of the eigenvalues of $p$ independent matrices.
For $p\ge n$, this reduces to the earlier result, since then we have
\[
U(n)^p \sim \bigoplus_{0\le i<n} U(1) \oplus \bigoplus_{n\le i<p} U(0),
\]
the latter component being trivial.  Similarly, the eigenvalue distribution
of a power of a random orthogonal or symplectic matrix can also be
described as a union of eigenvalue distributions.

After showing these results, we then consider the case of more general
compact Lie groups, especially since the above independence result extends
to this case.  It turns out in general that the eigenvalue density after
raising to a power can be expressed in terms of a certain parabolic
subgroup of an associated affine Weyl group.  (See \cite{Bourbaki},
\cite{Humphreys}, and \cite{MacdonaldIG:1972} for definitions and results
for affine Weyl groups.)  As a special case, we find that for a simple
compact Lie group $G$, the eigenvalues of $G^p$ are (essentially)
independent whenever $p$ is greater than the Coxeter number of $G$,
thus refining the result of \cite{Rains:1997a}.

\section{The classical groups}

We first consider the general problem of, given an eigenvalue distribution,
determining the distribution of powers of the eigenvalues.  We let
$\lambda_1$ through $\lambda_n$ denote complex numbers of norm 1, and
define
\[
dT := \prod_{0\le k<n} \frac{d\lambda_k}{2\pi i \lambda_k},
\]
so that $dT$ is the uniform density on the unit torus $T$.
The following lemma is straightforward:

\begin{lem}\label{lem:1.1}
Let $\pi(\lambda)$ be a Laurent polynomial in the variables $\lambda_j$
such that
\[
\pi(\lambda) dT
\]
is a probability density on the unit torus.  Then for any (monic) Laurent
monomial $\mu(\lambda)$,
\[
\int_T \pi(\lambda) \overline{\mu(\lambda)} dT
\]
is equal to the coefficient of $\mu(\lambda)$ in $\pi(\lambda)$.
\end{lem}

The following ``Main Lemma'' is crucial to our later results:

\begin{lem}
Fix an integer $p\ge 1$.
With $\pi(\lambda)$ as above, the joint density of the random variables
$\lambda^p_j$ is given by
\[
\pi^{(p)}(\lambda^{1/p}) dT,
\]
where $\pi^{(p)}(\lambda)$ is the sum of the $p$-divisible monomials of
$\pi(\lambda)$; that is, monomials in which each exponent is a multiple of
$p$.
\end{lem}

\begin{proof}
Since the torus is a compact set, it suffices to show that all joint
moments agree.  In other words, we need to show that if $\mu(\lambda)$
is a monic Laurent monomial, then
\[
\int_T \mu(\lambda) \pi^{(p)}(\lambda^{1/p}) dT
=
\int_T \mu(\lambda^p) \pi(\lambda) dT.
\]
But this follows immediately from Lemma \ref{lem:1.1}.
\end{proof}

We now obtain the power-map theorem for the unitary group:

\begin{thm}
For any integers $n,p\ge 1$,
\[
U(n)^p \sim \bigoplus_{0\le i<p} U(\lceil \frac{n-i}{p}\rceil);
\]
that is to say that the following two eigenvalue
distributions are the same:
\[
\eig(U^p)\sim \cup_{0\le i<p} \eig(U_i),
\]
where $U$ is a uniform random element of $U(n)$, and for $0\le i<p$
$U_i$ is an indepedent uniform random element of
$U(\lceil \frac{n-i}{p}\rceil)$.
\end{thm}

\begin{proof}
Recall \cite{Weyl} that the eigenvalue distribution of $U(n)$ is
given by the density
\[
{1\over n!} \sum_{\pi_1,\pi_2\in S_n}
\sigma(\pi_2\pi_1^{-1})
\prod_{0\le k<n} \lambda_k^{\pi_1(k)-\pi_2(k)}
dT.
\]
Equivalently, if we view $S_n$ as acting on the eigenvalues as
\[
\pi\cdot \lambda_i = \lambda_{\pi^{-1}(i)},
\]
we have the density
\[
{1\over n!} \sum_{\pi'\in S_n}
\pi'\cdot
(
\sum_{\pi\in S_n}
\sigma(\pi)
\prod_{0\le k<n} \lambda_k^{k-\pi(k)}
).
\]

By the Main Lemma, we must extract the $p$-divisible monomials.
Clearly, the action of $\pi'$ preserves divisibility, so
we may restrict our attention to
\[
q(\lambda)
=
\sum_{\pi\in S_n}
\sigma(\pi)
\prod_{0\le k<n} \lambda_k^{k-\pi(k)}.
\]
A monomial here is $p$-divisible if and only if
\[
\pi(k)\equiv k\pmod{p},\ \fa 0\le k<n.
\]
Given such a permutation, and given a congruence class $i\bmod p$,
we can define a new permutation $\pi^{(i)}\in S_{\lceil
\frac{n-i}{p}\rceil}$ by
\[
\pi^{(i)}(k) = \frac{\pi(pk+i)-i}{p}.
\]
Then the permutation $\pi$ is determined by the permutations $\pi^{(i)}$,
and we have
\begin{align}
\sigma(\pi) &= \prod_{0\le i<p} \sigma(\pi^{(i)})\\
\prod_{0\le k<n} \lambda_k^{k-\pi(k)}
&=
\prod_{0\le i<p}
\prod_{0\le k<\lceil\frac{n-i}{p}\rceil}
\lambda_{pk+i}^{pk-p\pi^{(i)}(k)}.
\end{align}
It follows that $q^{(p)}(\lambda^{1/p})$ factors:
\[
q^{(p)}(\lambda^{1/p})
=
\prod_{0\le i<p}
\left(
\sum_{\pi\in S_{\lceil\frac{n-i}{p}\rceil}}
\sigma(\pi)
\prod_{0\le k<\lceil\frac{n-i}{p}\rceil}
\lambda_{pk+i}^{k-\pi(k)}
\right).
\]
Since each factor has the form associated to the density of the unitary
group of dimension $\lceil\frac{n-i}{p}\rceil$, the result follows.
\end{proof}

As a special case, we recover a result of \cite{Rains:1997a}:

\begin{cor}
Fix an integer $n\ge 1$, and let $U$ be a random element of $U(n)$.
For any integer $p\ge n$, the eigenvalues of $U^p$ are i.i.d. uniform.
\end{cor}

\begin{proof}
Indeed, by the theorem,
\[
U(n)^p \sim \bigoplus_{0\le i<p} U(\lceil \frac{n-i}{p}\rceil)
=
\bigoplus_{0\le i<n} U(1)
\oplus
\bigoplus_{n\le i<p} U(0).
\]
But this is precisely the desired result.
\end{proof}

For the orthogonal group, the situation is somewhat more complicated,
both because the orthogonal group is not connected, and because
random orthogonal matrices are sometimes forced to have eigenvalues
$\pm 1$.  There are essentially four cases, depending on whether the
dimension is even or odd, and on whether the determinant is $1$ or $-1$.

\begin{thm}
Fix integers $n\ge 1$, $p\ge 1$.  Then we have the following identities
of eigenvalue distributions.
For $p$ odd:
\begin{align}
O^{\pm}(2n)^p
&\sim
O^{\pm}(2n_0)
\oplus
\bigoplus_{0\le i<(p-1)/2}
\Re U(\lceil\frac{2(n-n_0-i)}{p-1}\rceil).\\
O^{\pm}(2n+1)^p
&\sim
\bigoplus_{0\le i<(p-1)/2}
\Re U(\lceil\frac{2(n-n_1-i)}{p-1}\rceil)
\oplus
O^{\pm}(2n_1+1).
\end{align}
For $p$ even:
\begin{align}
O^{\pm}(2n)^p
&\sim
O^{\pm}(2n_0)
\oplus
\bigoplus_{0\le i<(p-2)/2}
\Re U(\lceil\frac{2(n-n_0-n_1-i)}{p-2}\rceil)
\oplus
O^{\mp}(2n_1+1).\\
O^{\pm}(2n+1)^p
&\sim
\bigoplus_{0\le i<p/2}
\Re U(\lceil\frac{2(n-i)}{p}\rceil).
\end{align}
Here $n_0=\lceil\frac{n}{p}\rceil$, $n_1=\lceil\frac{n-\lfloor
p/2\rfloor}{p}\rceil$, $O^\pm(n)$ represents the coset of $O(n)$ with
determinant $\pm 1$, $\Re U(n)$ represents the image of $U(n)$ in the
natural representation in $O(2n)$, and eigenvalues $\pm 1$ should be
ignored.
\end{thm}

\begin{proof}
Again, referring to \cite{Weyl}, the eigenvalue density for a random
orthogonal matrix is (up to an overall constant, and ignoring fixed
eigenvalues) given by
\[
\sum_{\rho'\in B_n}
\rho'\cdot
\left(
\sum_{\rho\in B_n}
\sigma(\rho)
\lambda^{\delta-\rho\cdot\delta}
\right)
dT
\]
where $B_n$ is the hyperoctahedral group (signed permutations),
$\sigma$ is a certain character of $B_n$, $\delta$ is a vector in
$\Z[1/2]^n$, acted on in the obvious way by $B_n$, and we define
\[
\lambda^v = \prod_{0\le k<n} \lambda_k^{v_k}.
\]
The ingredients $\sigma$ and $\delta$ are determined as follows:
\begin{align}
O^+(2n)\!:\,&\delta=(0,1,2\dots n-1),\ \sigma=\sigma_1\\
O^-(2n+2)\!:\,&\delta=(1,2,3\dots n),\ \sigma=\sigma_2\\
O^+(2n+1)\!:\,&\delta=(1/2,3/2,\dots n-1/2),\ \sigma=\sigma_2\\
O^-(2n+1)\!:\,&\delta=(1/2,3/2,\dots n-1/2),\ \sigma=\sigma_1
\end{align}
where $\sigma_1$ is the composition of the sign
character of $S_n$ with the natural projection $B_n\to S_n$,
and $\sigma_2$ is the natural sign character of $B_n$.

Consider a signed permutation $\rho$ with $\rho\cdot \delta\equiv
\delta\pmod{p}$.  Consider $\rho$ as a permutation of
\[
S = \{-n,1-n,\dots -2,-1,1,2,\dots n-1,n\},
\]
and define $\delta(x)$ for $x\in S$ by $\delta(k)=\delta_k$,
$\delta(-k)=-\delta_k$.  (Thus, for instance, for $O^+(2n)$, $\delta(k) =
|k|-1$.)  Then $\rho$ must preserve the partition of $S$ induced by
$\delta(x)\bmod p$.  The action of $\rho$ on an individual piece of the
partition is either as $S_n$, if $\delta(x)\not\equiv \delta(-x)\bmod p$,
or as $B_n$ if $\delta(x)\equiv
\delta(-x)\bmod p$.  Thus, the $p$-divisible part of
the inner sum
\[
\sum_{\rho\in B_n}
\sigma(\rho)
\lambda^{\delta-\rho\cdot\delta}
\]
factors; we obtain one $S_n$-type factor for each pair $\{i\bmod p,-i\bmod
p\}$ with $i\ne -i\bmod p$, and one $B_n$-type factor for each $i\bmod p$
with $i\equiv -i\bmod p$.  It remains only to determine the nature of the
factors (in particular, the character of $B_n$ that occurs) we obtain in
the eight cases.  This is an easy case-by-case analysis which we omit.

We thus deduce the following identities.
For $p$ odd:
\begin{align}
O^{\pm}(2n)^p
&\sim
O^{\pm}(2\lceil\frac{n}{p}\rceil)
\oplus
\bigoplus_{1\le i\le (p-1)/2}
\Re U(\lceil\frac{n-i}{p}\rceil + \lceil\frac{n+i}{p}\rceil-1).\\
O^{\pm}(2n+1)^p
&\sim
\bigoplus_{0\le i<(p-1)/2}
\Re U(\lceil\frac{n-i}{p}\rceil + \lceil\frac{n+1+i}{p}\rceil-1)
\oplus
O^{\pm}(2\lceil\frac{n-(p-1)/2}{p}\rceil+1).
\end{align}
For $p$ even:
\begin{align}
O^{\pm}(2n)^p
&\sim
O^{\pm}(2\lceil\frac{n}{p}\rceil)
\oplus
\bigoplus_{1\le i<p/2}
\Re U(\lceil\frac{n-i}{p}\rceil + \lceil\frac{n+i}{p}\rceil-1)
\oplus
O^{\mp}(2\lceil\frac{n-p/2}{p}\rceil+1).\\
O^{\pm}(2n+1)^p
&\sim
\bigoplus_{0\le i<p/2}
\Re U(\lceil\frac{n-i}{p}\rceil + \lceil\frac{n+1+i}{p}\rceil-1).
\end{align}

The theorem then follows by simple $\lceil\rceil$-manipulation.
\end{proof}

\begin{rem}
Since
\[
Sp(2n)\sim O^-(2n+2),
\]
ignoring the $\pm 1$ eigenvalues, this result also tells us the images of
$Sp(2n)$ under power maps.
\end{rem}

\begin{cor}
Fix an integer $n\ge 1$, and let $U$ be a random element of $O^\pm(n)$.
For any integer $p\ge n-1$, the eigenvalues other than $\pm 1$ of $U^p$ are
i.i.d. uniform conjugate pairs.  Similarly, if $U$ is a random element of
$Sp(2n)$ and $p\ge 2n+1$, then the eigenvalues of $U^p$ are i.i.d. uniform
conjugate pairs.
\end{cor}

\begin{proof}
It suffices to observe that the conclusion is true whenever
$O^\pm(n)^p$ is equivalent to a union of cosets $O^\pm(m)$ for $m\le 2$
and $\Re(U(m))$ for $m\le 1$.  By a case-by-case analysis, this holds
whenever $p\ge n-1$.
\end{proof}

Given the appearance of $\Re(U(n))$ above, it is appropriate to mention
the following relation:

\begin{thm}
For any $n\ge 0$,
\[
\Re(U(n)) \sim O^+(n+1)\oplus O^-(n+1),
\]
ignoring eigenvalues $\pm 1$.
\end{thm}

\begin{proof}
We claim that it suffices to prove
\[
g(1)g(-1)
\int_{U\in U(n)} \det(g(U))\det(g(\overline{U}))
=
\int_{O\in O^+(n+1)} \det(g(O))
\int_{O\in O^-(n+1)} \det(g(O)),
\label{eq:U=OO1}
\]
where $g(z)$ is an arbitrary function on the unit circle.  Indeed,
we may take $g$ to have the form
\[
g(z) = \prod_{1\le i\le m} (1-x_i z),
\]
at which point comparing coefficients of the $x_i$ tells us that
all joint moments of the polynomials
$(1-\lambda^2)\det(1-\lambda U)$
and
$\det(1-\lambda O_1)\det(1-\lambda O_2)$
agree, where $U\in U(n)$, $O_1\in O^+(n+1)$, $O_2\in O^-(n+1)$ are
uniform and independent.  But this implies the desired result.

To prove \eqref{eq:U=OO1}, one can express the integrals as determinants
(see, e.g., Theorems 2.1 and 2.2 of \cite{Rains:1999}) then use the
main lemma of \cite{WilfHS:1992}.  An alternate proof, given as
Corollary 2.4 of \cite{Rains:1999}, involves expressing the integrals
in terms of orthogonal polynomials on the unit circle.
\end{proof}

\begin{rem}
It is possible to give similar proofs of the the other results of this
paper; in particular, one uses the fact that the orthogonal polynomials
with respect to a weight $g(z^p)$ are simply related to the orthogonal
polynomials with respect to $g(z)$.  However, it is unclear how
to apply this approach to nonclassical groups, just as it is unclear how
to apply the other approach to prove this result.
\end{rem}

\section{Congruential subgroups of Weyl groups}

For the unitary and orthogonal groups, the key observation was that
an appropriate subgroup of the Weyl group turned out to be itself a
product of Weyl groups.  More precisely, we were given a Weyl group $W$,
a lattice $\Lambda$, and a vector $v$ in the ambient space of $\Lambda$,
and considered the group
\[
W^{v+\Lambda}:=\{\rho:\rho\in W|\rho(v)-v\in\Lambda\}.
\]
For instance, for $O^+(2n+1)$, we had
\begin{align}
v &= (1/2,3/2,\dots (n-1)/2),\\
\Lambda &= p \Z^n.
\end{align}
Equivalently, we could divide $v$ and $\Lambda$ by $p$, thus making
$\Lambda$ equal to the root lattice of $O^+(2n+1)$.  This suggests
that we should first study $W^{v+\Lambda}$ when $\Lambda$ is the root
lattice of $W$.

\begin{thm}
Let $W$ be a (finite) Weyl group (not necessarily simple), and let
$\Lambda_a$ be the root lattice of $W$.  Then for any vector $v\in\R\Lambda_a$,
$W^{v+\Lambda_a}$ is the Weyl group generated by the roots of $W$ it contains.
\end{thm}

\begin{proof}
The key observation is that an element $\rho\in W$ is in $W^{v+\Lambda_a}$
if and only if there exists a translation $t_\lambda\in \Lambda_a$ such that
\[
t_\lambda^{-1}(\rho(v)) = v.
\]
Thus instead of considering the stabilizer in $W$ of the translate
$v+\Lambda$, we can consider the stabilizer in $W^+:=W\ltimes\Lambda_a$
of the vector $v$; we have a canonical isomorphism between $W^{v+\Lambda_a}$
and $(W^+)^v$.

Since $\Lambda_a$ is the root lattice of $W$, $W^+$ is an affine Weyl group.
But then we can apply proposition V.3.3.2 of \cite{Bourbaki} to conclude
that $(W^+)^v$ is generated by the reflections of $W^+$ that fix $v$.
The theorem follows immediately.
\end{proof}

\begin{cor}
With $W$, $v$, $\Lambda_a$ and $W^+$ as above, $W^{v+\Lambda_a}$ is isomorphic
to a finite parabolic subgroup of $W^+$.
\end{cor}

\begin{proof}
Indeed, $(W^+)^v$ is a finite subgroup of $W^+$ generated by reflections,
so by definition is parabolic.
\end{proof}

It will be helpful to refine the above result somewhat.  Fix a fundamental
chamber of $W$.  This then forces a choice of fundamental chamber of $W^+$,
as well as choices of fundamental chambers for the corresponding parabolic
subgroups.  We recall the notation of \cite{MacdonaldIG:1972}: given an
element $\rho\in W^+$, $D\rho$ is the element $t_{-\rho(0)} \rho$ which,
since it fixes $0$, is in $W$.

\begin{lem}
There exists an element $\rho\in W^+$ such that $\rho(v)$ is in
the fundamental chamber of $W^+$ and $(D\rho)(v)$ is in the fundamental
chamber of $W^{\rho(v)+\Lambda_a}\cong W^{v+\Lambda_a}$.  The resulting
points $\rho(v)$ and $(D\rho)(v)$ are then independent of $\rho$.
\end{lem}

\begin{proof}
Certainly, there exists an element $\rho_0\in W^+$ such that $\rho_0(v)$ is
in the fundamental chamber of $W^+$; the point $\rho_0(v)$ is then
independent of the choice of $\rho_0$.  The remaining freedom in the choice
of $\rho_0$ is simply that we can apply any element of $(W^+)^{\rho_0(v)}$.
Since
\[
D((W^+)^{\rho_0(v)}) = W^{\rho_0(v)+\Lambda}
\]
there exists an element $\rho_1$ of $(W^+)^{\rho_0(v)}$ such that
\[
(D\rho_1)(D\rho_0)(v)
\]
is in the fundamental chamber of $W^{\rho_0(v)+\Lambda_a}$, and again
the resulting point is unique.  Taking $\rho=\rho_1\rho_0$, we
are done.
\end{proof}

Given the point $v$, we define new points $\overline{v}$ and $\tilde{v}$
by
\[
\overline{v}:= \rho(v)\ \text{and}\ \tilde{v}:=(D\rho)(v),
\]
with $\rho$ as above.  Note that $\tilde{v}-\overline{v}\in\Lambda_a$.

If $\Lambda$ is not the root lattice, then $W\ltimes \Lambda$ is in general
no longer an affine Weyl group, so the above results do not apply.
We say that $\Lambda$ is a ``subweight'' lattice of $W$ if it contains the
roots of $W$, and satisfies
\[
\frac{2\langle \Lambda,r\rangle}{|r|^2} \subset \Z
\]
for all roots $r$.  We recall:

\begin{lem}
Let $\Lambda$ be a subweight lattice of $W$ contained in $\R\Lambda_a$.
Then for each coset of $(W\ltimes \Lambda)/W^+$, there exists a unique
representative that preserves the fundamental chamber of $W^+$, giving an
injection from $\Lambda/\Lambda_a$ into the group of automorphisms of the
Dynkin diagram of $W^+$.
\end{lem}

This gives the following result:

\begin{thm}
Let $\Lambda$ be a subweight lattice of $W$, and define $\Lambda_0 :=
\Lambda\cap (\R\Lambda_a)$.  Then for any vector $v\in \R^n$,
\[
W^{v+\Lambda_a} \unlhd W^{v+\Lambda};
\]
the quotient is isomorphic to the subgroup of $\Lambda_0/\Lambda_a$ such
that the corresponding transformations of the fundamental chamber preserve
$\overline{v}$.
\end{thm}

\begin{proof}
The main complication is the fact that $\R\Lambda_a$ might not equal
$\R\Lambda$.  If not, let $V$ be the orthogonal complement of $\R\Lambda_a$
in $\R\Lambda$, and write $v=v_0+v_1$ with $v_0\in \R\Lambda_a$ and $v_1\in
V$.  Then we can write
\[
\rho(v)-v = \rho(v_0)-v_0\in \R\Lambda_a.
\]
Thus $\rho(v)-v\in\Lambda$ if and only if $\rho(v)-v\in\Lambda_0$.

We may therefore assume that $\Lambda\subset\R\Lambda_a$, and thus
$\Lambda_0=\Lambda$.  Furthermore, since the desired result is invariant
under conjugation by $W^+$, we may assume that $v=\overline{v}$.
Fix a coset of $W^+$ in $W\ltimes \Lambda$, and let $\psi$ be the
representative of that coset that preserves the fundamental chamber.
If $\psi(v)=v$, then $\psi$ clearly normalizes $(W^+)^v$, and
\[
(W\ltimes\Lambda)^v \cap \psi W^+ = \psi (W^+)^v.
\]
On the other hand, if $\psi(v)\ne v$ then
\[
(W\times \Lambda)^v \cap \psi W^+ = 0.
\]
The result follows.
\end{proof}

Thus given $\overline{v}$, the group $W^{\overline{v}+\Lambda}$ can
essentially be read off by inspection.

\section{General compact Lie groups}

Let $G$ be a connected compact Lie group, choose a maximal torus $T$ of
$G$, and let $W$ be the corresponding Weyl group.  Choosing an isomorphism
$T\cong U(1)^n$, we obtain a set of $n$ characters $\lambda_i$ of $T$ such
that for any representation $R$ of $G$, the eigenvalues of a matrix in
$R(T)$ are given by appropriate (monic) Laurent monomials in the
$\lambda_i$.  Since every element of $G$ can be conjugated into $T$, the
characters $\lambda_i$ can be thought of as ``eigenvalue generators''.

There is a one-to-one correspondence between Laurent monomials in the
eigenvalue generators and vectors in $\Lambda^G:=\Z^n$, given by
\[
v\mapsto \lambda^v := \prod_{1\le i\le n} \lambda_i^{v_i}.
\]
The Killing form, Weyl group (in its action on $T$), and root lattice of
$G$ thus all induce corresponding structures in $\Lambda^G$.  We find
that $\Lambda^G$ is a subweight lattice of $W$.

Before considering the eigenvalue distribution of $G$, it will be helpful
to digress briefly on the relation between $G$ and $\Lambda^G$.  Let
$\tilde{\Lambda}$ denote the lattice formed by adjoining all weights of $W$
to $\Lambda^G$, and let $\Lambda_2$ and $\Lambda_1$ denote the projections
of $\Lambda$ onto the root space of $W$ and its orthogonal complement,
respectively.  Then $\Lambda_1\times\Lambda_2$ is again a subweight lattice
of $W$, and in fact decomposes as a product of integer lattices (with null
inner product) and simple weight lattices.  Thus $\Lambda_1\times \Lambda_2
=
\Lambda^{G^+}$ for a group $G^+$ of the form
\[
U(1)^n \times \tilde{H}
\]
where $\tilde{H}$ is simply connected.
We can then write $G$ as the quotient of $G^+$ by some discrete subgroup
of its center.  Indeed, this subgroup is given as a subgroup of $T(G^+)$
as the kernel of the natural projection
\[
T(G^+) = (\Lambda^{G^+})^*\to (\Lambda^G)^*=T(G),
\]
and is thus naturally isomorphic to $(\Lambda^{G^+}/\Lambda^G)^*$.

Now, the density of the eigenvalue distribution of $G$ (with respect to
Haar measure) has the form $f(\lambda) dT$ where
\[
f(\lambda)
\propto
\sum_{\rho'\in W}
\rho'\cdot
\left(
\sum_{\rho\in W}
\sigma(\rho)
\lambda^{\delta-\rho\delta}
\right),
\]
where $\delta$ is the Weyl vector of $W$, i.e. half the sum of the
positive roots.  We thus find that
\begin{align}
f^{(p)}(\lambda^{1/p})
&\propto
\sum_{\rho'\in W}
\rho'\cdot
\left(
\sum_{\rho\in W^{\delta/p+\Lambda^G}}
\sigma(\rho)
\lambda^{(\delta-\rho\delta)/p}
\right)\\
&=
\sum_{\rho'\in W}
\rho'\cdot
\left(
\sum_{\rho\in W^{(p)}}
\sigma(\rho)
\lambda^{(\tilde{\delta}-\rho\tilde{\delta})/p}
\right).\\
&\propto
\sum_{\rho'\in W}
\rho'\cdot
\left(
\sum_{\rho_1,\rho_2\in W^{(p)}}
\sigma(\rho_1\rho_2^{-1})
\lambda^{(\rho_1\tilde{\delta}-\rho_2\tilde{\delta})/p}
\right),
\end{align}
where we define
\[
W^{(p)} := W^{\tilde{\delta}/p+\Lambda^G}.
\]

The inner sum thus looks roughly like the eigenvalue density of some
different Lie group.  For it to actually {\it be} an eigenvalue density (or
rather, for our methods to {\it prove} it an eigenvalue density\footnote{In
particular, we observe that $O^+(2n)^p$ and $Sp(2n)^p$ violate these
conditions!}), we need three things to happen.  First, we need
$\Lambda^G$ to be a subweight lattice of $W^{(p)}$; luckily, this is
trivial, since every root of $W^{(p)}$ is a root of $W$.  Second, we need
\[
W^{(p)}=W^{\tilde{\delta}/p+\Lambda_a},
\]
so that the stabilizer is actually a Weyl group.  Finally, we need the
projection of $\tilde{\delta}/p$ to the root space of
$W^{(p)}$ to be equal to the Weyl vector of $W^{(p)}$.

These last two conditions are, unfortunately, not always satisfied,
although we can at least readily determine when they are; see
for instance Section \ref{sec:tables} below.
One partial result is:

\begin{thm}
The center of $G$ is connected if and only if $\Lambda_0=\Lambda_a$.  Thus
if the center of $G$ is connected (in particular if $G$ is adjoint), then
\[
W^{(p)}=W^{\tilde{\delta}/p+\Lambda_a},
\]
and $W^{(p)}$ is a Weyl group.
\end{thm}

\begin{proof}
Writing $G=G^+/Z$ as above, we observe that $G$ has connected center if
and only if for every element $z\in Z(\tilde{H})$, there exists some
element of $Z$ that projects to $z$.  Dualizing, this is precisely the
requirement that $\Lambda_0=\Lambda_a$.
\end{proof}

Thus, at least in the connected center case, the only obstacle is
the third condition.  That this is, indeed, a problem can be seen from
the tables of Section \ref{sec:tables}; this is the main obstacle to
a truly satisfying result in general.

Similar remarks hold if $G$ is not connected.  Indeed, much of the
structure theory can be extended to this case (see \cite{Rains:1997a}
and \cite{dSiebenthalJ:1956}\footnote{The author was unaware of the results in
\cite{dSiebenthalJ:1956} when writing \cite{Rains:1997a}, and thus most of the
structural results of \cite{Rains:1997a} had already appeared (with different
proofs) in \cite{dSiebenthalJ:1956}, with the notable exceptions of the density
formula and the independence result (see below).}).  A connected component
of a compact Lie group is specified by a pair $(G_0,a)$, where $G_0$ is a
connected compact Lie group and $a$ is an automorphism of $G_0$ (up to
conjugation by outer automorphisms and multiplication by inner
automorphisms; in particular, $a$ may be chosen to have finite order).
Using the classification of simple Lie groups, one readily obtains
the following irreducible local possibilities, each indexed by a positive
integer $n$:

\begin{itemize}
\item{${}^nU$:} $G_0=U(1)^{\phi(n)}$, $a$ satisfies the cyclotomic
polynomial of order $n$.
\item{${}^nH$:} $G_0=H^n$, with $H$ simple; $a$ acts as a cyclic shift.
\item{${}^nA^{(2)}_m$:} $G_0=(SU(m+1))^n$; $a$ acts as a cyclic shift, twisted
by the outer automorphism of $SU(m+1)$.
\item{${}^nD^{(2)}_m$:} $G_0=(\tilde{SO}(2m))^n$; $a$ acts as a cyclic shift,
twisted by the (classical) outer automorphism of $\tilde{SO}(2m)$.
\item{${}^nE^{(2)}_6$:} $G_0=(E_6)^n$; $a$ acts as a cyclic shift, twisted by
the outer automorphism of $E_6$.
\item{${}^nD^{(3)}_4$:} $G_0=(\tilde{SO}(8))^n$; $a$ acts as a cyclic shift,
twisted by the triality automorphism of $\tilde{SO}(8)$.
\end{itemize}

The case ${}^nU$ is essentially trivial; either $n=1$, in which case the
(single) eigenvalue generator is uniformly distributed, or $n>1$, in which
case there are no eigenvalue generators (because the eigenvalues are
constant over the component).  In the other cases, if we let $f_X(\lambda)$
denote the density for ${}^1X$, then the density for ${}^nX$ is given by
$f_X(\lambda^n)$.  This implies:

\begin{cor}
Let $X$ either denote a simply connected, compact, simple Lie group or one
of $A^{(2)}_m$, $D^{(2)}_m$, $E^{(2)}_6$, or $D^{(3)}_4$.  Then for all
positive integers $n$ and $p$, the eigenvalue density of $({}^n X)^p$
is given by
\[
f(\lambda^{n/\gcd(n,p)}),
\]
where $f$ is the density for $({}^1 X)^{p/\gcd(n,p)}$.
\end{cor}

Thus it suffices to consider the five ``interesting'' cases with $n=1$.  In
each case, we find that the density has the expected form corresponding to
an appropriate Weyl group (keeping the same Weyl vector), but using a
different weight lattice in place of the root lattice.  The effective Weyl
group and $\Lambda^G$ are given as follows:
\begin{itemize}
\item{$A^{(2)}_{2m-1}$:} $W=B_m$.  $\Lambda^G=(1/2) D_m$ when
$G_0$ has even fundamental group, and $(1/2) Z^m$ when $G_0$ has odd
fundamental group.
\item{$A^{(2)}_{2m}$:} $W=C_m$.  $\Lambda^G=(1/2) Z^m$.
\item{$D^{(2)}_m$:} $W=C_{m-1}$.  $\Lambda^G=Z^{m-1}$ for the adjoint and
orthogonal groups, and $\Lambda^G=C^\perp_{m-1}$ otherwise.
\item{$E^{(2)}_6$:} $W=F_4$. $\Lambda^G = F_4^\perp$.
\item{$D^{(3)}_4$:} $W=G_2$. $\Lambda^G = G_2^\perp$.
\end{itemize}

In each case, if we define $\Lambda_a$ to be the weight lattice
corresponding to the adjoint group, we find that $W\ltimes \Lambda_a$ is an
affine Weyl group, and the above results carry over.  We observe the
following relations in the adjoint case:
\[
{}^n A^{(2)}_{2m} \sim {}^{2n} D^{(2)}_{m+1} \sim {}^{2n} C_m.
\]

\section{Independence results}

In \cite{Rains:1997a}, it was shown that for any connected component $C$ of any
compact Lie group $G$, there exists a threshold $P$ such that for any $p>P$
and any representation of $G$, we have the relation
\[
C^p \sim U(1)^r
\]
on the eigenvalue generators; that is, the eigenvalue generators of $C^p$
are independent and uniform.  Using the above considerations, we can now
give an explicit value to the threshold:

\begin{thm}
Let $G$ be a connected compact Lie group, let $h$ be the maximum Coxeter
number of the simple factors of $W(G)$, and let $r$ be the rank of $G$.  Then
\[
p>h   \implies G^p\sim U(1)^r
\]
and conversely
\[
G^p\sim U(1)^r \implies p\ge h.
\]
If the center of $G$ is connected, then $C^h\sim U(1)^r$.
\end{thm}

\begin{proof}
The key observation is that $\delta$ is not stabilized by any element of
$W$.  Consequently, the sum
\[
\sum_{\rho\in W^{(p)}} \sigma(\rho) \lambda^{\delta-\rho(\delta)}
\]
is equal to 1 precisely when the group $W^{(p)}$ is trivial.  If the
sum is 1, then clearly $G^p\sim U(1)^r$; conversely, if $G^p\sim U(1)^r$,
then the sum must equal 1.

It thus remains to consider the group $W^{(p)}$.  Aside from diagram
automorphisms, $W^{(p)}$ is the product of the groups corresponding to the
simple factors of $G$; we may thus assume that $G$ is simple.  In this
case, we can explicitly verify that when $p<h$, $W^{(p)}$ is nonempty (it
is straightforward for the classical groups, and a short computation for
the exceptional groups).  For $p\ge h$, we have the following lemma (from
Prop. 7.3 of \cite{MacdonaldIG:1972}):

\begin{lem}
Let $W$ be a simple (finite) Weyl group, and let $W^+=W\ltimes \Lambda_a$.
If $h$ is the Coxeter number and $\delta$ the Weyl vector of $W$, then
$\delta/h$ is the centroid of the fundamental chamber of $W^+$, and is thus
invariant under all diagram automorphisms of $W^+$.
\end{lem}

Thus when $p\ge h$, $\delta/p$ is strictly in the interior of $W^+$, so
any nontrivial element of $W^{(p)}$ must come from a diagram automorphism;
since this is impossible when the center of $G$ is connected, we obtain the
desired result for $p=h$.  When $p>h$, $\delta/p$ treats the highest root
of $W^+$ differently from the other roots, so any diagram automorphism
preserving $\delta/p$ must preserve the highest root.  But this precludes
the diagram automorphisms corresponding to the cosets of the root lattice
in the weight lattice, giving the desired result.
\end{proof}

A similar result holds in the disconnected case (with an appropriate
definition of Coxeter number), with $p$ replaced by $p/\gcd(p,n)$, except
for ${}^nA^{(2)}_m$ and ${}^nE^{(2)}_6$, when $p$ must be replaced by
$p/\gcd(p,2n)$, and for ${}^nD^{(3)}_4$, when $p$ must be replaced by
$p/\gcd(p,3n)$.

Thus for instance, we find that
\[
E_8(n)^p \sim U(1)^8
\]
precisely when $p\ge 30$.
Similarly,
\[
SU(n)^p\sim U(1)^{n-1}
\]
precisely when $p>n$, since $SU(n)$ is not adjoint (and we readily verify
that $W^{(n)}$ has $n$ elements in this case); this threshold was
incorrectly given as $n-1$ in \cite{Rains:1997a}.

\section{Tables}\label{sec:tables}

We conclude the paper by giving a table of $\wt{\delta/p}$ and
$\overline{\delta/p}$ (for $1\le p\le h$) for the exceptional groups; from
this, it is straightforward to read off the appropriate power-relations.

Vectors in the root space are specified by their inner products with the
fundamental affine roots, scaled by the length of the root, except that we use
$\overline{k}$ to represent $1/p-k$.  We, in fact, give only
$\wt{\delta/p}$ via this scheme; $\overline{\delta/p}$ is then obtained by
replacing $k\mapsto 0$, $\overline{k}\mapsto \overline{0}$ (and is thus
essentially given by the locations of the bars).  The Dynkin
diagram of $W^{(p)}$ is then read off as the subdiagram spanned by the
indices without bars, and the projection of $\wt{\delta/p}$ to the root
space of $W^{(p)}$ is a Weyl vector precisely when all of those indices are
1.  Finally, we append an asterisk when $\overline{\delta/p}$ is
symmetric under the appropriate diagram transformation.

We remark that it is not at all clear why this encoding scheme should
happen to work!

\break
%
% my_sym:
%
%    If the argument is a positive number, typeset $\mathbf{#1}$.
%    If the argument is a nonnegative number, typeset $\overline{#1}$.
%
%  This will be used in the preambles to the haligns in this section.
%
\def\my_sym#1{\relax\ifnum#1<1 $\overline{\number-#1}$\else%
\ifnum#1<99 $\mathbf{\number#1}$\else \hphantom{$\overline{99}$}\fi\fi}
  
\begin{multicols}{2}
\begin{itemize}
\item{$G_2$}
% 0-1<<<2
$
\begin{smallmatrix}
0-1{<}\equiv 2
\end{smallmatrix}
$
\halign{\hfil$#$&\hfil$#$&&\hfil\my_sym{#}\hfil&$\,#$\hfil\cr
\noalign{\vskip-\baselineskip}
&&          99&&99&&99&&99&&99&&99&&99\cr
p={}& 1: (&-6&& 1&& 1 & )\cr
p={}& 2: (& 2&&-3&& 1 & )\cr
p={}& 3: (& 0&&-2&& 1 & )\cr
p={}& 4: (& 0&& 1&&-1 & )\cr
p={}& 5: (& 1&&-1&& 0 & )\cr
p={}& 6: (&-1&& 0&& 0 & )\cr
}
\item{$F_4$}
% 0-1-2<<3-4
$
\begin{smallmatrix}
0-1-2\Leftarrow 3-4
\end{smallmatrix}
$
\halign{\hfil$#$&\hfil$#$&&\hfil\my_sym{#}\hfil&$\,#$\hfil\cr
\noalign{\vskip-\baselineskip}
&&          99&&99&&99&&99&&99&&99&&99\cr
p={}& 1: (&-12&& 1&& 1&& 1&& 1 & )\cr
p={}& 2: (& 2&&-6&& 1&& 1&& 1 & )\cr
p={}& 3: (& 1&& 2&&-4&& 1&& 1 & )\cr
p={}& 4: (& 0&& 1&&-3&& 1&& 1 & )\cr
p={}& 5: (& 1&& 1&&-2&& 1&&-1 & )\cr
p={}& 6: (&-1&& 1&&-2&& 1&& 0 & )\cr
p={}& 7: (& 0&&-1&& 1&&-1&& 1 & )\cr
p={}& 8: (& 0&& 0&&-1&& 1&&-1 & )\cr
p={}& 9: (& 0&& 0&& 1&&-1&& 0 & )\cr
p={}&10: (& 0&& 1&&-1&& 0&& 0 & )\cr
p={}&11: (& 1&&-1&& 0&& 0&& 0 & )\cr
p={}&12: (&-1&& 0&& 0&& 0&& 0 & )\cr
}
\item{$E_6$}
% 0-1-6;2-3-6;4-5-6
$
\begin{smallmatrix}
&&0&&\\
&&\mid&&\\
&&1&&\\
&&\mid&&\\
2-{}&3-{}&6&{}-5&{}-4
\end{smallmatrix}
$
\medskip\medskip
\halign{\hfil$#$&\hfil$#$&&\hfil\my_sym{#}\hfil&$\,#$\hfil\cr
\noalign{\vskip-\baselineskip}
&&          99&&99&&99&&99&&99&&99&&99\cr
p={}& 1{:}\, (&-12&& 1&& 1&& 1&& 1&& 1&& 1&)  \cr
p={}& 2{:}\, (&  2&&-6&& 1&& 1&& 1&& 1&& 1&)  \cr
p={}& 3{:}\, (&  1&& 2&& 1&& 1&& 1&& 1&&-4&)^*\cr
p={}& 4{:}\, (&  0&& 1&& 1&& 1&& 1&& 1&&-3&)  \cr
p={}& 5{:}\, (&  1&& 1&&-1&& 1&&-1&& 1&&-2&)  \cr
p={}& 6{:}\, (& -1&& 1&& 0&& 1&& 0&& 1&&-2&)^*\cr
p={}& 7{:}\, (&  0&&-1&& 1&&-1&& 1&&-1&& 1&)  \cr
p={}& 8{:}\, (&  0&& 0&&-1&& 1&&-1&& 1&&-1&)  \cr
p={}& 9{:}\, (&  0&& 0&& 0&&-1&& 0&&-1&& 1&)^*\cr
p={}&10{:}\, (&  0&& 1&& 0&& 0&& 0&& 0&&-1&)  \cr
p={}&11{:}\, (&  1&&-1&& 0&& 0&& 0&& 0&& 0&)  \cr
p={}&12{:}\, (& -1&& 0&& 0&& 0&& 0&& 0&& 0&)^*\cr
}
\item{$E_7$}
% 0-1-2-3-4;3-5-6-7
$
\begin{smallmatrix}
&&&4\\
&&&\mid\\
0-{}&1-{}&2-{}&3&{}-5&{}-6&{}-7\\
\end{smallmatrix}
$
\smallskip
\halign{\hfil$#$&\hfil$#$&&\hfil\my_sym{#}\hfil&$\,#$\hfil\cr
\noalign{\vskip-\baselineskip}
&&          99&&99&&99&&99&&99&&99&&99\cr
p={}& 1: (& 1&& 1&& 1&& 1&& 1&& 1&& 1&&-18 & )\cr
p={}& 2: (& 1&& 1&& 1&& 1&&-9&& 1&& 1&& 2 & )^*\cr
p={}& 3: (& 1&& 1&& 1&& 1&& 1&&-6&& 2&& 1 & )\cr
p={}& 4: (&-3&& 1&& 1&& 1&& 1&&-4&& 1&& 1 & )\cr
p={}& 5: (& 1&& 1&& 1&&-3&& 1&& 1&& 1&&-2 & )\cr
p={}& 6: (&-1&& 1&& 1&&-3&& 1&& 1&& 1&& 0 & )^*\cr
p={}& 7: (& 1&&-1&& 1&&-2&& 1&& 1&& 1&&-2 & )\cr
p={}& 8: (&-1&& 1&& 1&&-2&& 1&& 1&&-1&& 0 & )\cr
p={}& 9: (& 1&&-1&& 1&&-2&& 1&& 1&& 0&& 0 & )\cr
p={}&10: (&-1&& 1&&-1&& 1&&-1&&-1&& 1&& 0 & )^*\cr
p={}&11: (& 0&&-1&& 1&&-1&& 0&& 1&&-1&& 1 & )\cr
p={}&12: (& 0&& 0&&-1&& 1&& 0&&-1&& 1&&-1 & )\cr
p={}&13: (& 0&& 0&& 0&&-1&& 1&& 1&&-1&& 0 & )\cr
p={}&14: (& 0&& 0&& 0&& 1&&-1&&-1&& 0&& 0 & )^*\cr
p={}&15: (& 0&& 0&& 1&&-1&& 0&& 0&& 0&& 0 & )\cr
p={}&16: (& 0&& 1&&-1&& 0&& 0&& 0&& 0&& 0 & )\cr
p={}&17: (& 1&&-1&& 0&& 0&& 0&& 0&& 0&& 0 & )\cr
p={}&18: (&-1&& 0&& 0&& 0&& 0&& 0&& 0&& 0 & )^*\cr
}
\item{$E_8$}
% 0-1-2-3-4-5-6;5-7-8
$
\begin{smallmatrix}
&&&&&6\\
&&&&&\mid\\
0-{}&1-{}&2-{}&3-{}&4-{}&5&{}-7&{}-8\\
\end{smallmatrix}
$
\smallskip
\halign{\hfil$#$&\hfil$#$&&\hfil\my_sym{#}\hfil&$\,#$\hfil\cr
\noalign{\vskip-\baselineskip}
&&          99&&99&&99&&99&&99&&99&&99\cr
p={}& 1: (&-30&& 1&& 1&& 1&& 1&& 1&& 1&& 1&& 1 & )\cr
p={}& 2: (& 2&& 1&& 1&& 1&& 1&& 1&& 1&& 1&&-15 & )\cr
p={}& 3: (& 1&& 1&& 1&& 1&& 1&& 1&&-10&& 1&& 2 & )\cr
p={}& 4: (& 2&& 1&& 1&&-7&& 1&& 1&& 1&& 1&& 1 & )\cr
p={}& 5: (& 1&& 1&& 1&& 2&&-6&& 1&& 1&& 1&& 1 & )\cr
p={}& 6: (& 0&& 1&& 1&& 1&&-5&& 1&& 1&& 1&& 1 & )\cr
p={}& 7: (& 1&&-2&& 1&& 1&&-4&& 1&& 1&& 1&& 1 & )\cr
p={}& 8: (& 2&&-3&& 1&& 1&& 1&&-3&& 1&& 1&& 1 & )\cr
p={}& 9: (& 0&&-2&& 1&& 1&& 1&&-3&& 1&& 1&& 1 & )\cr
p={}&10: (& 0&& 1&&-1&& 1&& 1&&-3&& 1&& 1&& 1 & )\cr
p={}&11: (& 1&& 1&&-2&& 1&& 1&&-2&& 1&& 1&&-1 & )\cr
p={}&12: (&-1&& 1&&-2&& 1&& 1&&-2&& 1&& 1&& 0 & )\cr
p={}&13: (& 0&&-1&& 1&&-1&& 1&&-2&& 1&& 1&& 1 & )\cr
p={}&14: (& 0&& 0&&-1&& 1&& 1&&-2&& 1&& 1&&-1 & )\cr
p={}&15: (& 0&& 0&& 1&&-1&& 1&&-2&& 1&& 1&& 0 & )\cr
p={}&16: (& 0&& 1&&-1&& 1&&-1&& 1&&-1&&-1&& 1 & )\cr
p={}&17: (& 1&&-1&& 1&&-1&& 1&&-1&& 0&& 1&&-1 & )\cr
p={}&18: (&-1&& 1&&-1&& 1&&-1&& 1&& 0&&-1&& 0 & )\cr
p={}&19: (& 0&&-1&& 1&&-1&& 1&&-1&& 1&& 0&& 0 & )\cr
p={}&20: (& 0&& 0&&-1&& 1&&-1&& 1&&-1&& 0&& 0 & )\cr
p={}&21: (& 0&& 0&& 0&&-1&& 1&&-1&& 0&& 1&& 0 & )\cr
p={}&22: (& 0&& 0&& 0&& 0&&-1&& 1&& 0&&-1&& 1 & )\cr
p={}&23: (& 0&& 0&& 0&& 0&& 0&&-1&& 1&& 1&&-1 & )\cr
p={}&24: (& 0&& 0&& 0&& 0&& 0&& 1&&-1&&-1&& 0 & )\cr
p={}&25: (& 0&& 0&& 0&& 0&& 1&&-1&& 0&& 0&& 0 & )\cr
p={}&26: (& 0&& 0&& 0&& 1&&-1&& 0&& 0&& 0&& 0 & )\cr
p={}&27: (& 0&& 0&& 1&&-1&& 0&& 0&& 0&& 0&& 0 & )\cr
p={}&28: (& 0&& 1&&-1&& 0&& 0&& 0&& 0&& 0&& 0 & )\cr
p={}&29: (& 1&&-1&& 0&& 0&& 0&& 0&& 0&& 0&& 0 & )\cr
p={}&30: (&-1&& 0&& 0&& 0&& 0&& 0&& 0&& 0&& 0 & )\cr
}
\item{$(D^{(3)}_4)^3$}
% 0-1>>>2
$
\begin{smallmatrix}
0-1\equiv {>}2
\end{smallmatrix}
$
\halign{\hfil$#$&\hfil$#$&&\hfil\my_sym{#}\hfil&$\,#$\hfil\cr
\noalign{\vskip-\baselineskip}
&&          99&&99&&99&&99&&99&&99&&99\cr
p={}& 1: (& 1&& 2 &&-6& )\cr
p={}& 2: (& 0&& 1 &&-3& )\cr
p={}& 3: (& 1&&-1 && 0& )\cr
p={}& 4: (&-1&& 0 && 0& )\cr
}
\item{$(E^{(2)}_6)^2$}
% 0-1-2>>3-4
$
\begin{smallmatrix}
0-1-2\Rightarrow 3-4
\end{smallmatrix}
$
\halign{\hfil$#$&\hfil$#$&&\hfil\my_sym{#}\hfil&$\,#$\hfil\cr
\noalign{\vskip-\baselineskip}
&&          99&&99&&99&&99&&99&&99&&99\cr
p={}& 1: (& 2&& 1&& 1&& 2&&-12 & )\cr
p={}& 2: (& 1&& 1&& 2&&-6&&  2 & )\cr
p={}& 3: (& 0&& 1&& 1&&-4&&  2 & )\cr
p={}& 4: (&-2&& 1&& 1&&-2&&  0 & )\cr
p={}& 5: (& 0&&-1&& 1&&-2&&  2 & )\cr
p={}& 6: (& 0&& 0&& 1&&-2&&  0 & )\cr
p={}& 7: (& 0&& 1&&-1&& 0&&  0 & )\cr
p={}& 8: (& 1&&-1&& 0&& 0&&  0 & )\cr
p={}& 9: (&-1&& 0&& 0&& 0&&  0 & )\cr
}
\end{itemize}
\end{multicols}

\end{document}